\documentclass[10pt,draft,twoside,leqno,a4paper]{article}
\usepackage{amsmath,amsthm,amsfonts,amssymb,verbatim}
\newtheorem{Theorem}{Theorem}
\newtheorem{Lemma}{Lemma}

\newtheorem{Proposition}{Proposition}

\newcommand{\dash}{\mathchoice
    {\mkern.70mu\raise.50ex\hbox{\vrule height.1ex width.40em depth0pt}}
    {\mkern.40mu\raise.48ex\hbox{\vrule height.1ex width.30em depth0pt}}
    {\mkern.33mu\raise.30ex\hbox{\vrule height.1ex width.25em depth0pt}}
    {\mkern.10mu\raise.20ex\hbox{\vrule height.1ex width.20em depth0pt}}
    }

\newcommand{\R}{\mathbb R}

\newcommand{\N}{\mathbb N}
\newcommand{\eps}{\varepsilon}

\newcommand{\weak}{\rightharpoondown}
\newcommand{\bra}{\langle}
\newcommand{\ket}{\rangle}
\renewcommand{\d}{{\mathrm d}}
\begin{document}
\title{Mountain pass solutions for
a mean field equation from two-dimensional turbulence}
\author{Tonia Ricciardi\thanks{
Supported in part by 
the MIUR National Project {\em Variational Methods and
Nonlinear Differential Equations}.}
\\
{\small Dipartimento di Matematica e Applicazioni
``R.~Caccioppoli"}\\
{\small Universit\`a di Napoli Federico II}\\
{\small Via Cintia, 80126 Napoli, Italy}\\
{\small fax: +39 081 675665}\\
{\small\tt{ tonia.ricciardi@unina.it}}\\
}
\date{December 5, 2006}
\maketitle
\begin{abstract}
Using Struwe's ``monotonicity trick"
and the recent blow-up analysis of Ohtsuka and Suzuki, 
we prove the existence of mountain pass solutions to a
mean field equation arising in two-dimensional turbulence.
\end{abstract}
\begin{description}
\item {\textsc{Key Words:}} mean field equation, mountain pass solution
\item {\textsc{MSC 2000 Subject Classification:}} 35J60
\end{description}
\section{Introduction and main result}
We consider the following problem:
\begin{align}
\label{pb}
&-\Delta_gu=\lambda_1\left(\frac{e^u}{\int_Me^u\,\d v_g}-\frac{1}{|M|}\right)
-\lambda_2\left(\frac{e^{-u}}{\int_Me^{-u}\,\d v_g}-\frac{1}{|M|}\right)\\
\nonumber
&\int_M u\,\d v_g=0,
\end{align}
where $(M,g)$ is a compact, orientable, Riemannian 2-manifold without boundary,
$\Delta_g$ denotes the Laplace-Beltrami operator,
$|M|$ denotes the volume of $M$
and where $\lambda_1,\lambda_2$ are positive constants.
Equation~\eqref{pb} arises in the context of the statistical mechanics description of 
two-dimensional turbulence,
see Joyce and Montgomery~\cite{JM} and Pointin and Lundgren~\cite{PL}.
Problem~\eqref{pb} has a variational structure. Indeed, solutions
to \eqref{pb} correspond to critical points for the functional
\begin{align*}
\label{funct}
I_{\lambda_1,\lambda_2}&(u)=\\
=&\frac{1}{2}\int_M|\nabla_gu|^2\,\d v_g
-\lambda_1\ln\left(\frac{1}{|M|}\int_M e^u\,\d v_g\right)
-\lambda_2\ln\left(\frac{1}{|M|}\int_M e^{-u}\,\d v_g\right)
\end{align*}
defined on the Sobolev space
\[
E=\{u\in H^1(M)\ :\ \int_Mu\,\d v_g=0\},
\]
equipped with the norm $\|u\|^2=\int_M|\nabla_gu|^2\d v_g$.
A detailed analysis of the functional $I_{\lambda_1,\lambda_2}$ 
has been recently carried out by Ohtsuka and Suzuki~\cite{OS}, together with a blow-up analysis for
problem~\eqref{pb}. In particular, it is shown in \cite{OS} that 
$I_{\lambda_1,\lambda_2}$ is bounded below 
if and only if
$0\le\lambda_1,\lambda_2\le8\pi$, thus improving the classical sharp Moser-Trudinger inequality.
The corresponding problem under Dirichlet boundary
conditions has been recently considered in Sawada, Suzuki and Takahashi~\cite{SST}.
\par
We note that \eqref{pb} always admits the trivial solution $u\equiv0$. 
Our aim in this note is to 
prove the existence of nontrivial solutions to problem~\eqref{pb}
for suitable values of $\lambda_1,\lambda_2$.
We denote by $\mu_1(M)$ the first nonzero eigenvalue of $\Delta_g$, namely,
\[
\mu_1(M)=\inf_{u\in E\setminus\{0\}}\frac{\int_M|\nabla_g u|^2\,\d v_g}{\int_Mu^2\,\d v_g}.
\]
We assume that $\mu_1(M)$ satisfies the condition
\begin{equation}
\label{eigenvalueassumption}
8\pi<\mu_1(M)|M|<16\pi.
\end{equation}
Condition~\eqref{eigenvalueassumption} is satisfied, e.g., on the flat torus~$\R^2/\mathbb Z^2$,
since $\mu_1(\R^2/\mathbb Z^2)=4\pi^2$.
We define
\[
\Lambda=\left\{(\lambda_1,\lambda_2)\in\R_0^+\times\R_0^+\,:\,
\lambda_1+\lambda_2<\mu_1(M)|M|\ \mathrm{and}\ \max\{\lambda_1,\lambda_2\}>8\pi\right\}.
\]
Note that condition~\eqref{eigenvalueassumption} ensures that $\Lambda\neq\emptyset$.
In fact, $\Lambda$ is the union of the triangle with vertices $(8\pi,0)$, $(\mu_1(M)|M|,0),
(8\pi,\mu_1(M)|M|-8\pi)$, and of its 
reflection with respect to the straight line $\lambda_1=\lambda_2$.
We prove the following result:
\begin{Theorem}
\label{thm:main}
Suppose $(M,g)$ satisfies condition~\eqref{eigenvalueassumption}.
Then problem~\eqref{pb} admits a nonzero solution for every $(\lambda_1,\lambda_2)\in\Lambda$.
\end{Theorem}
We note that when $\lambda_2=0$, Theorem~\ref{thm:main} 
reduces to Theorem~2.1 obtained by Struwe and Tarantello in \cite{ST}.
In fact, the proof of Theorem~\ref{thm:main}
employs the variational approach introduced in \cite{ST}
and further developed by Jeanjean in \cite{Je}, together with
a fine blow-up analysis for solutions to \eqref{pb} by Ohtsuka and Suzuki~\cite{OS}. 
More precisely, in Section~\ref{sec:prop} we show that for every $(\lambda_1,\lambda_2)\in\Lambda$
the functional $I_{\lambda_1,\lambda_2}$ has a mountain pass structure. However, since
$\max\{\lambda_1,\lambda_2\}>8\pi$ we have a lack of compactness, and general Palais-Smale 
sequences may ``blow-up" in the sense of Brezis and Merle~\cite{BM}. 
Therefore, in Section~\ref{sec:proof} 
we employ Struwe's ``monotonicity trick"~\cite{St} in order to find a 
bounded Palais-Smale sequence 
for almost every $(\lambda_1,\lambda_2)\in\Lambda$, and consequently a solution to \eqref{pb}.
Finally, we exploit the blow-up analysis of Ohtsuka and Suzuki~\cite{OS}
to obtain solutions for every $(\lambda_1,\lambda_2)\in\Lambda$.
Related results in the context of $SU(3)$~Toda systems may be found in Chae, Ohtsuka and Suzuki~\cite{COS}
and in Lucia and Nolasco~\cite{LN}.
\section{Some properties of $I_{\lambda_1,\lambda_2}$}
\label{sec:prop}
In this section we prove some properties of $I_{\lambda_1,\lambda_2}$
which will be useful in what follows.
In particular, we prove the following result:
\begin{Proposition}
\label{prop:mp}
The functional $I_{\lambda_1,\lambda_2}$
has a mountain pass structure for every $(\lambda_1,\lambda_2)\in\Lambda$.
\end{Proposition}
Before proving Proposition~\ref{prop:mp} we introduce some notation and
we establish some preliminary results.
For every $u\in E$ we set
\begin{equation}
\label{Gdef}
G(u)=\ln\left(\frac{1}{|M|}\int_M e^u\,\d v_g\right).
\end{equation}
Then, we may write 
\[
I_{\lambda_1,\lambda_2}(u)=\frac{1}{2}\int_M|\nabla_g u|^2\,\d v_g
-\lambda_1 G(u)-\lambda_2G(-u).
\]
For every $\phi,\psi\in E$ we have
\begin{align*}
G'(u)\phi=&\frac{\int_M e^u\phi\,\d v_g}{\int_M e^u\,\d v_g}\\
\bra G''(u)\phi,\psi\ket
=&\frac{\int_M e^u\phi\psi\,\d v_g\int_Me^u\,\d v_g-\int_M e^u\phi\,\d v_g\int_M e^u\psi\,\d v_g}
{\left(\int_M e^u\,\d v_g\right)^2}.
\end{align*}
In particular, recalling that $\int_M\phi\,\d v_g=0$, we obtain
\begin{equation}
\label{Gexp}
G'(0)=0,\qquad\qquad\bra G''(0)\phi,\phi\ket=\frac{\int_M\phi^2\,\d v_g}{|M|}.
\end{equation}
Furthermore, by the Cauchy-Schwarz inequality, we have
\[
\left(\int_M e^u\phi\,\d v_g\right)^2\le\int_M e^u\phi^2\,\d v_g\int_M e^u\,\d v_g
\]
for every $u,\phi\in E$, and therefore
\begin{equation}
\label{CS}
\bra G''(u)\phi,\phi\ket\ge0\qquad\qquad\forall u,\phi\in E.
\end{equation}
We shall also use the following property of $G$:
\begin{Lemma}
\label{lem:G}
There exists a constant $C_M>0$ depending on $M$ only, such that
\[
\|G'(u)\|\le C_M\exp\left\{\frac{1}{8\pi}\|u\|^2\right\}
\]
for all $u\in E$.
\end{Lemma}
\begin{proof}
By Jensen's inequality we have $\int_Me^u\d v_g\ge|M|$.
Therefore, for every $\phi\in E$ we have:
\[
|G'(u)\phi|\le\frac{1}{|M|}\int_Me^u\phi\d v_g
\le\frac{1}{|M|}\left(\int_Me^{2u}\d v_g\right)^{1/2}\left(\int_M\phi^2\d v_g\right)^{1/2}.
\]
In view of the Moser-Trudinger inequality:
\begin{equation}
\label{mt}
\int e^u\d v_g\le C_0\exp\left\{\frac{1}{16\pi}\|u\|^2\right\}
\qquad\forall u\in E,
\end{equation}
where $C_0>0$ depends on $M$ only,
and of the Poicar\'e inequality
\begin{equation}
\label{poincare}
\int_M\phi^2\,\d v_g\le\frac{1}{\mu_1(M)}\int_M|\nabla_g\phi|^2\,\d v_g
\qquad\forall\phi\in E,
\end{equation}
we derive
\[
|G'(u)\phi|\le\frac{\left[C_0\exp\left\{\frac{1}{4\pi}\|u\|^2\right\}\right]^{1/2}\|\phi\|}
{\mu_1(M)^{1/2}|M|}
\]
and the asserted estimate follows with $C_M=\left(C_0/\mu_1(M)\right)^{1/2}|M|^{-1}$. 
\end{proof}
\begin{Lemma}
\label{lem:locmin}
For every $\lambda_1,\lambda_2\ge0$ such that $\lambda_1+\lambda_2<\mu_1(M)|M|$
the function $u\equiv0$ is a local minimum for $I_{\lambda_1,\lambda_2}$.
\end{Lemma}
\begin{proof}
We have:
\begin{align*}
I_{\lambda_1,\lambda_2}(0)=&0,\qquad\qquad I_{\lambda_1,\lambda_2}'(0)=0,\\
\bra I_{\lambda_1,\lambda_2}''(0)\phi,\phi\ket=&\int_M|\nabla_g\phi|^2\,\d v_g
-\frac{\lambda_1+\lambda_2}{|M|}\int_M\phi^2\,\d v_g.
\end{align*}
In view of the Poicar\'e inequality~\eqref{poincare}
we derive
\begin{align*}
\bra I_{\lambda_1,\lambda_2}''(0)\phi,\phi\ket\ge\left(1-\frac{\lambda_1+\lambda_2}{\mu_1(M)|M|}\right)\|\phi\|^2.
\end{align*}
Now the claim follows by Taylor expansion of $I_{\lambda_1,\lambda_2}$ at $u\equiv0$.
\end{proof}
Let $p_0\in M$ and let $r_0>0$ be a constant smaller than the injectivity
radius of $M$ at $p_0$. Let $\mathcal B_{r_0}=\{p\in M:d_g(p,p_0)<r_0\}$
denote the geodesic ball of radius $r_0$ centered at $p_0$.
For every $\eps>0$ let $v_\eps$ be the function defined by
\[
v_\eps(p)=\begin{cases}
\ln\frac{\eps^2}{(\eps^2+d_g(p,p_0)^2)^2}&\mathrm{in\ }\mathcal B_{r_0}\\
\ln\frac{\eps^2}{(\eps^2+r_0^2)^2}&\mathrm{in\ }M\setminus \mathcal B_{r_0}
\end{cases}
\]
and let $u_\eps\in E$ be correspondingly defined by 
\[
u_\eps=v_\eps-\frac{1}{|M|}\int_M v_\eps\,\d v_g.
\]
We have the following asymptotic expansions.
\begin{Lemma}
\label{lem:ue}
There exists $r_0>0$ such that the following asymptotic expansions hold, 
with respect to $\eps\to0$:
\begin{align}
\label{nablaue}
&\int_M|\nabla_g u_\eps|^2\,\d v_g=32\pi\ln\frac{1}{\eps}+O(1)\\
\label{eue}
&\ln\int_M e^{u_\eps}\,\d v_g=2\ln\frac{1}{\eps}+O(1)\\
\label{e-ue}
&\ln\int_M e^{-u_\eps}\,\d v_g=O(1).
\end{align}
\end{Lemma}
\begin{proof}
In geodesic coordinates centered at $p$ we have 
$g_{ij}(x)=\delta_{ij}+O(|x|^2)$ for every $x\in B_{r_0}=\{x\in\R^2:|x|<r_0\}$. 
Consequently, $\d v_g=(1+O(|x|^2))\,\d x$ and
$|\nabla_g u|^2\,\d v_g=|\nabla u|^2(1+O(|x|^2))\,\d x$, 
where $\nabla$ denotes differentiation with respect to the Euclidean metric.
Moreover, identifying $v_\eps,u_\eps$ with their pullbacks to $B_{r_0}$
under the geodesic coordinate system,
we have in $B_{r_0}$:
\begin{align*}
v_\eps(x)=&\ln\frac{\eps^2}{(\eps^2+|x|^2)^2},\\
\nabla u_\eps(x)=&\nabla v_\eps(x)=-\frac{4x}{\eps^2+|x|^2}.
\end{align*}
Proof of \eqref{nablaue}. We have
\begin{align*}
\int_{B_{r_0}}|\nabla u_\eps|^2\,\d x=&32\pi\int_0^{r_0}\frac{r^3}{(\eps^2+r^2)^2}\,\d r\\
=&32\pi\left(\int_0^{r_0}\frac{r\,\d r}{\eps^2+r^2}-\eps^2\int_0^{r_0}\frac{r\,\d r}{(\eps^2+r^2)^2}\right)\\
=&16\pi\left(\left.\ln(\eps^2+r^2)\right|_0^{r_0}+\left.\frac{\eps^2}{\eps^2+r^2}\right|_0^{r_0}\right)
=16\pi\ln\eps^{-2}+O(1).
\end{align*}
Moreover, 
\begin{align*}
&\int_{B_{r_0}}|\nabla u_\eps|^2|x|^2\,\d x\le 16\int_{B_{r_0}}\d x=O(1).
\end{align*}
Since $\nabla_gu_\eps\equiv0$ in $M\setminus\mathcal B_{r_0}$, 
we conclude that
\[
\int_M|\nabla_gu_\eps|^2\,\d v_g=\int_{B_{r_0}}|\nabla u_\eps|^2(1+O(|x|^2))\,\d x
=16\pi\ln\eps^{-2}+O(1)
\]
and \eqref{nablaue} is established.
\par
Proof of \eqref{eue}.
We have:
\begin{align*}
&\int_{B_{r_0}}e^{v_\eps}\,\d x=\int_{B_{r_0}}\frac{\eps^2}{(\eps^2+|x|^2)^2}\,\d x
=\int_{B_{r_0/\eps}}\frac{\d y}{(1+|y|^2)^2}=\int_{\R^2}\frac{\d y}{(1+|y|^2)^2}+\circ(1)\\
&\int_{B_{r_0}}e^{v_\eps}|x|^2\,\d x=\int_{B_{r_0}}\frac{\eps^2|x|^2}{(\eps^2+|x|^2)^2}\,\d x
=\eps^2\int_{B_{r_0/\eps}}\frac{|y|^2\,\d y}{(1+|y|^2)^2}=\circ(1),\\
&\int_{M\setminus\mathcal B_{r_0}}e^{v_\eps}\,\d v_g=\frac{\eps^2}{(\eps^2+r_0^2)}|M\setminus\mathcal B_{r_0}|=\circ(1).
\end{align*}
Therefore, we obtain
\begin{equation*}
\label{eve}
\int_Me^{v_\eps}\,\d v_g=\int_{B_{r_0}}e^{v_\eps}(1+O(|x|^2))\,\d x
+\int_{M\setminus\mathcal B_{r_0}}e^{v_\eps}\,\d v_g=\int_{\R^2}\frac{\d y}{(1+|y|^2)^2}+\circ(1)
\end{equation*}
and consequently
\begin{equation}
\label{lneve}
\ln\int_Me^{v_\eps}\,\d v_g=\ln\left(\int_{\R^2}\frac{\d y}{(1+|y|^2)^2}+\circ(1)\right)=O(1).
\end{equation}
Moreover, we have
\begin{align*}
\int_Mv_\eps\,\d v_g=&|M|\ln\eps^2-2\int_{B_{r_0}}\ln(\eps^2+|x|^2)\,\d v_g
-2|M\setminus\mathcal B_{r_0}|\ln(\eps^2+r_0^2)\\
=&|M|\ln\eps^2+O(1)
\end{align*}
and therefore the mean value of $v_\eps$ satisfies
\begin{equation}
\label{intve}
\frac{1}{|M|}\int_M v_\eps\,\d v_g=\ln\eps^2+O(1).
\end{equation}
Consequently, in view of \eqref{lneve} and \eqref{intve},
\begin{align*}
\ln\int_Me^{u_\eps}\,\d v_g=\ln\int_Me^{v_\eps}\,\d v_g-\frac{1}{|M|}\int_Mv_\eps\,\d v_g
=\ln\eps^{-2}+O(1).
\end{align*}
Hence, \eqref{eue} is established.
\par
Proof of \eqref{e-ue}.
Using the elementary formula 
$\int_{B_{r_0}}|x|^p\,\d x=2\pi(p+2)^{-1}r_0^{p+2}$, 
we compute:
\begin{align*}
\int_{B_{r_0}}e^{-v_\eps}\,\d x=&\int_{B_{r_0}}\frac{(\eps^2+|x|^2)^2}{\eps^2}\,\d x\\
=&\eps^{-2}\int_{B_{r_0}}|x|^4\,\d x+2\int_{B_{r_0}}|x|^2\,\d x+\eps^2\int_{B_{r_0}}\,\d x\\
=&\frac{\pi r_0^6}{3}\eps^{-2}+\pi r_0^4+\eps^2\pi r_0^2,\\
\int_{B_{r_0}}e^{-v_\eps}|x|^2\,\d x=&\eps^{-2}\int_{B_{r_0}}|x|^6\,\d x+2\int_{B_{r_0}}|x|^4
+\eps^2\int_{B_{r_0}}|x|^2\,\d x\\
=&\frac{\pi r_0^8}{4}\eps^{-2}+\frac{2\pi r_0^6}{3}+\eps^2\frac{\pi r_0^4}{2}.
\end{align*}
Therefore, 
\begin{align*}
\int_{B_{r_0}}e^{-v_\eps}\,\d v_g=&\int_{B_{r_0}}\frac{(\eps^2+|x|^2)^2}{\eps^2}(1+O(|x|^2))\,\d x\\
=&\frac{\pi r_0^6}{3}(1+O(r_0)^2)\eps^{-2}+\pi r_0^4(1+O(r_0)^2)
+\eps^2\pi r_0^2(1+O(r_0)^2)
\end{align*}
and by choosing $r_0>0$ sufficiently small, we obtain
\[
\int_{B_{r_0}}e^{-v_\eps}\,\d v_g=\beta_0\eps^{-2}+\beta_1+\circ(1)
\]
for some $\beta_0,\beta_1>0$.
On the other hand, in $M\setminus\mathcal B_{r_0}$ we have
\begin{align*}
\int_{M\setminus\mathcal B_{r_0}}e^{-v_\eps}\,\d v_g=&\frac{(\eps^2+r_0^2)^2}{\eps^2}|M\setminus\mathcal B_{r_0}|\\
=&|M\setminus\mathcal B_{r_0}|\eps^{-2}+2r_0^2|M\setminus\mathcal B_{r_0}|+\circ(1).
\end{align*}
It follows that
\begin{equation*}
\label{e-ve}
\int_Me^{-v_\eps}\,\d v_g=(\beta_0+|M\setminus\mathcal B_{r_0}|)\eps^{-2}+\beta_1
+2r_0^2|M\setminus\mathcal B_{r_0}|+\circ(1)
\end{equation*}
and consequently,
\begin{align*}
\label{lne-ve}
\ln\int_Me^{-v_\eps}\,\d v_g=&\ln\left[(\beta_0+|M\setminus\mathcal B_{r_0}|)\eps^{-2}
+\beta_1+2r_0^2|M\setminus\mathcal B_{r_0}|+\circ(1)\right]\\
\nonumber
=&\ln\eps^{-2}+\ln\left[\beta_0+|M\setminus\mathcal B_{r_0}|+\circ(1)\right]
=\ln\eps^{-2}+O(1).
\end{align*}
In view of the above and \eqref{intve}, we derive
\begin{align*}
\ln\int_Me^{-u_\eps}\,\d v_g=&\ln\int_Me^{-v_\eps}\,\d v_g+\frac{1}{|M|}\int_Mv_\eps\,\d v_g\\
=&\ln\eps^{-2}+\ln\eps^2+O(1)=O(1).
\end{align*}
This establishes \eqref{e-ue}.
\end{proof}
\begin{proof}[Proof of Proposition~\ref{prop:mp}]
Let $(\lambda_1,\lambda_2)\in\Lambda$.
In view of Lemma~\ref{lem:locmin}, $u\equiv0$ is a strict nondegenerate
local minimum for
$I_{\lambda_1,\lambda_2}$.
Moreover,
in view of the expansions \eqref{nablaue}, \eqref{eue} and \eqref{e-ue},
we have
\[
I_{\lambda_1,\lambda_2}(u_\eps)=2(8\pi-\lambda_1)\ln\frac{1}{\eps}+O(1).
\]
By symmetry,
\[
I_{\lambda_1,\lambda_2}(-u_\eps)=2(8\pi-\lambda_2)\ln\frac{1}{\eps}+O(1).
\]
In particular, it follows that whenever $\max\{\lambda_1,\lambda_2\}>8\pi$, the functional
$I_{\lambda_1,\lambda_2}$ is unbounded below. Consequently,
for all $(\lambda_1,\lambda_2)\in\Lambda$, $I_{\lambda_1,\lambda_2}$ has a mountain pass structure.
\end{proof}
\textit{Remark.} The unboundedness of $I_{\lambda_1,\lambda_2}$ 
when $\max\{\lambda_1,\lambda_2\}>8\pi$ was derived by Ohtsuka and Suzuki~\cite{OS}
from certain properties of a functional related to the Moser-Trudinger inequality~\eqref{mt}.
However, we have proved it directly,
since we shall use some properties of the function~$u_\eps$ in the next section.
\section{Proof of Theorem~\ref{thm:main}}
\label{sec:proof}
Using Struwe's ``monotonicity trick"~\cite{St} together with the blow-up analysis
for \eqref{pb} developed by Ohtsuka and Suzuki~\cite{OS},
we first prove the following result.
\begin{Proposition}
\label{prop:mono}
For a.e.\ $(\lambda_1,\lambda_2)\in\Lambda$ 
there exists a nonzero critical point
for $I_{\lambda_1,\lambda_2}$.
\end{Proposition}
In order to prove Proposition~\ref{prop:mono}, we begin by setting a minimax argument.
In view of Proposition~\ref{prop:mp} and Lemma~\ref{lem:ue}, 
there exists $\bar u=\bar u(\lambda_1,\lambda_2)\in E$
satisfying $\|\bar u\|\ge 1$ and $I_{\lambda_1,\lambda_2}(\bar u)<0$.
We consider the set of paths
\[
\Gamma=\{\gamma\in C([0,1],E)\ :\ \gamma(0)=0,\ \gamma(1)=\bar u\}.
\]
Then, the value
\[
c(\lambda_1,\lambda_2)=\inf_{\gamma\in\Gamma}\max_{t\in[0,1]}I_{\lambda_1,\lambda_2}(\gamma(t))
\]
is finite. 
More precisely, the following estimate holds.
\begin{Lemma}
\label{lem:c}
For every $\omega>0$ there exists $\rho_\omega>0$ independent of $(\lambda_1,\lambda_2)$
such that
\[
c(\lambda_1,\lambda_2)\ge\frac{\rho_\omega^2}{2}\left(1-\frac{\lambda_1+\lambda_2}{\mu_1(M)|M|}-\omega\right).
\]
\end{Lemma}
\begin{proof}
By Taylor expansion of $G$ at 0, where $G$ is the functional defined in \eqref{Gdef},
in view of \eqref{Gexp} we have:
\[
G(u)=\frac{\int_Mu^2\,\d v_g}{|M|}+\circ(\|u\|^2).
\]
Hence, recalling that $\lambda_1+\lambda_2<\mu_1(M)|M|$ and the inequality~\eqref{poincare}, we may write
\begin{align*}
I_{\lambda_1,\lambda_2}(u)=&\frac{1}{2}\int_M|\nabla_gu|^2\,\d v_g
-\frac{\lambda_1+\lambda_2}{2|M|}\int_Mu^2\,\d v_g+(\lambda_1+\lambda_2)\circ(\|u\|^2)\\
\ge&\frac{1}{2}\left(1-\frac{\lambda_1+\lambda_2}{\mu_1(M)|M|}\right)\|u\|^2-\mu_1(M)|M|\circ(\|u\|^2).
\end{align*}
Hence, for every $\omega>0$ there exists $\rho_\omega\in(0,1)$ independent of $(\lambda_1,\lambda_2)$
such that
\[
I_{\lambda_1,\lambda_2}(u)
\ge\frac{1}{2}\left(1-\frac{\lambda_1+\lambda_2}{\mu_1(M)|M|}-\omega\right)\|u\|^2
\]
for every $\|u\|\le\rho_\omega$,
and in particular
\[
I_{\lambda_1,\lambda_2}(u)
\ge\frac{\rho_\omega^2}{2}\left(1-\frac{\lambda_1+\lambda_2}{\mu_1(M)|M|}-\omega\right)
\]
for all $u\in E$ satisfying $\|u\|=\rho_\omega$.
Since for every $\gamma\in\Gamma$ there exists $\bar t$ for which $\|\gamma(\bar t)\|=\rho_\omega$,
we have
\[
\max_{t\in[0,1]}I_{\lambda_1,\lambda_2}(\gamma(t))\ge I_{\lambda_1,\lambda_2}(\gamma(\bar t))
\ge\frac{\rho_\omega^2}{2}\left(1-\frac{\lambda_1+\lambda_2}{\mu_1(M)|M|}-\omega\right).
\]
Now the claim follows recalling the definition of $c(\lambda_1,\lambda_2)$.
\end{proof}
We note that for every fixed $(\lambda_1,\lambda_2)\in\Lambda\setminus\partial\Lambda$
the function $\lambda\mapsto c(\lambda_1+\lambda,\lambda_2+\lambda)$
is well-defined and monotone decreasing for all
$\lambda$ near 0. Therefore it is differentiable at a.e.\ $\lambda$.
Consequently, the function $\eps\mapsto c(\lambda_1+\eps,\lambda_2+\eps)$
is differentiable at $\eps=0$ for a.e.\ $(\lambda_1,\lambda_2)\in\Lambda\setminus\partial\Lambda$.
In what follows, such a $(\lambda_1,\lambda_2)\in\Lambda$ will be fixed,
and to simplify notation, we set:
\begin{align*}
I=&I_{\lambda_1,\lambda_2},\\
I_\eps=&I_{\lambda_1+\eps,\lambda_2+\eps}=I-\eps(G(u)+G(-u))\\
c=&c_{\lambda_1,\lambda_2}\\
c_\eps=&c(\lambda_1+\eps,\lambda_2+\eps)\\
c'=&\left.\frac{dc_\eps}{d\eps}\right|_{\eps=0},
\end{align*}
where the functional $G$ is defined \eqref{Gdef}.
We consider the set $X_\eps=X_\eps(\lambda_1,\lambda_2)$ defined by
\[
X_\eps=\left\{u\in E\ :\ \begin{matrix}u=\gamma(t)\ \mathrm{for\ some\ } \gamma\in\Gamma\ \mathrm{s.t.\ }
\max_{t\in[0,1]}I(\gamma(t))\le c+\eps\\
\mathrm{and\ }I_\eps(u)\ge c_\eps-\eps
\end{matrix}
\right\}.
\]
The following estimate holds:
\begin{Lemma}
\label{lem:monoest}
Let $(\lambda_1,\lambda_2)\in\Lambda$ be such that
$c'(\lambda_1,\lambda_2)$ exists. Then
there exist constants $\eta>0$, $C_1>0$ such that,
for every $\eps\in(0,\eta)$ and for every $u\in X_\eps$,
there holds:
\begin{enumerate}
\item[(i)]
$\|u\|\le C_1$
\item[(ii)]
$c-(|c'|+2)\eps\le I(u)\le c+\eps$.
\end{enumerate}
\end{Lemma}
\begin{proof}
Since $c'$ exists, we have
\[
c_\eps=c+c'\eps+\circ(\eps).
\]
Proof of (i). For every $u\in X_\eps$ we have, as $\eps\to0$:
\begin{align*}
0\le&G(u)+G(-u)
=\frac{I(u)-I_\eps(u)}{\eps}\le\frac{c+\eps-(c_\eps-\eps)}{\eps}\\
=&-c'+2+\circ(1)=|c'|+2+\circ(1)
\end{align*}
and therefore there exists $\eta>0$ such that
\[
G(u)+G(-u)\le|c'|+3\qquad\qquad\forall u\in X_\eps,\ \forall\eps\in(0,\eta).
\]
It follows that 
\begin{align*}
\|u\|^2=&2[I(u)+\lambda_1G(u)+\lambda_2G(-u)]\le2[c+\eps+\max\{\lambda_1,\lambda_2\}(G(u)+G(-u))]\\
\le&2[c+\eta+\mu_1(M)|M|(|c'|+3)]
\end{align*}
and therefore the asserted estimate~(i) holds with $C_1^2=2[c+\eta+\mu_1(M)|M|(|c'|+3)]$.
\par\noindent
Proof of (ii).
We have, recalling the monotonicity property of $I_\eps$:
\begin{align*}
I(u)\ge I_\eps(u)\ge c_\eps-\eps
=c+c'\eps+\circ(\eps)-\eps\ge c-(|c'|+2)\eps.
\end{align*}
Since $I(u)\le c+\eps$ by definition of $X_\eps$, (ii)
is also established.
\end{proof}
Now we show that if $c'(\lambda_1,\lambda_2)$ exists,
then $X_\eps$ necessarily contains a bounded Palais-Smale sequence.
\begin{Lemma}
\label{lem:PS}
Let $(\lambda_1,\lambda_2)\in\Lambda$ be such that $c'(\lambda_1,\lambda_2)$
exists. Then
there exists a sequence $(u_n)_{n\in\N}$, $u_n\in E$, such that $\|u_n\|\le C_1$,
$I(u_n)\to c$, $\|I'(u_n)\|\to0$.
\end{Lemma}
\begin{proof}
If not, there exists $\delta>0$ such that
for every $u\in E$ satisfying $\|u\|\le C_1$ and $|I(u)-c|<\delta$ 
we have $\|I'(u)\|\ge\delta$.
In order to derive a contradiction, we need a suitable deformation.
Let $\eps_n\to0$ and let $\varphi:\R\to\R$ be a smooth cutoff function satisfying
$0\le\varphi\le1$, $\varphi(s)=1$ for all $s\ge-1$, $\varphi(s)=0$ for all $s\le-2$,
and let
\[
\Phi_n(u)=\varphi\left(\frac{I_{\eps_n}(u)-c_{\eps_n}}{\eps_n}\right).
\]
Then $\Phi_n(u)=1$ for all $u$ such that $I_{\eps_n}(u)\ge c_{\eps_n}-\eps_n$ 
and in view of Lemma~\ref{lem:c} we have $\Phi_n(0)=\Phi_n(\bar u)=0$
for sufficiently large $n$.
Consequently, setting
\[
\widetilde\gamma(t)=\gamma(t)-\sqrt{\eps_n}\,\Phi_n(\gamma(t))
\frac{I_{\eps_n}'(\gamma(t))}{\|I_{\eps_n}'(\gamma(t))\|},
\]
we have $\widetilde\gamma\in\Gamma$ for every $\gamma\in\Gamma$. 
For every $n\in\N$, let $\gamma_n\in\Gamma$ be such that $\max_{t\in[0,1]}I(\gamma_n(t))\le c+\eps_n$.
We note that in view of Lemma~\ref{lem:monoest}, 
for all $u\in X_{\eps_n}$ with sufficiently large $n$ we have $\|u\|\le C_1$,
$|I(u)-c|<\delta$ and therefore, by assumption, we have $\|I'(u)\|\ge\delta$.
\par\noindent
CLAIM.
For sufficiently large $n$, we have
\begin{equation}
\label{max}
\max_{t\in[0,1]}I_{\eps_n}(\widetilde\gamma_n(t))
=\max_{\gamma_n(t)\in X_{\eps_n}}I_{\eps_n}(\widetilde\gamma_n(t)).
\end{equation}
Proof of Claim. By definition of $c_{\eps_n}$ we have
\begin{equation}
\label{maxequiv}
\max_{t\in[0,1]}I_{\eps_n}(\widetilde\gamma_n(t))
=\max_{\{t:I_{\eps_n}(\widetilde\gamma_n(t))\ge c_{\eps_n}-\eps_n/2\}}I_{\eps_n}(\widetilde\gamma_n(t)).
\end{equation}
We note that in view of \eqref{CS} we have, for every $u,v\in E$:
\begin{align*}
\bra I_{\eps_n}''(u)v,v\ket=&\int_M|\nabla_g v|^2\,\d v_g-(\lambda_1+\eps_n)\bra G''(u)v,v\ket
-(\lambda_2+\eps_n)\bra G''(-u)v,v\ket\\
\le&\int_M|\nabla_g v|^2\,\d v_g
\end{align*}
and therefore, an expansion to the second order of $I_{\eps_n}$ yields 
\begin{align}
\nonumber
I_{\eps_n}(u+v)
=&I_{\eps_n}(u)+I_{\eps_n}'(u)v+\frac{1}{2}\bra I_{\eps_n}''(u+sv)v,v\ket\\
\label{exp}
\le&I_{\eps_n}(u)+I_{\eps_n}'(u)v+\frac{1}{2}\|v\|^2,
\end{align}
where $s\in(0,1)$.
Using \eqref{exp}, for all $t\in[0,1]$ such that $I_{\eps_n}(\widetilde\gamma_n(t))>c_{\eps_n}-\eps_n/2$
we estimate:
\begin{align*}
c_{\eps_n}-\frac{\eps_n}{2}<&I_{\eps_n}(\widetilde\gamma_n(t))\\
\le&I_{\eps_n}(\gamma_n(t))+I_{\eps_n}'\,(\gamma_n(t))\left(\widetilde\gamma_n(t)-\gamma_n(t)\right)
+\frac{1}{2}\|\widetilde\gamma_n(t)-\gamma_n(t)\|^2\\
=&I_{\eps_n}(\gamma_n(t))-\sqrt{\eps_n}\Phi_n(\gamma_n(t))\|I_{\eps_n}'(\gamma_n(t))\|
+\frac{\eps_n}{2}\Phi_n^2(\gamma_n(t))\\
\le&I_{\eps_n}(\gamma_n(t))+\frac{\eps_n}{2}.
\end{align*}
It follows that
\[
I_{\eps_n}(\gamma_n(t))\ge c_{\eps_n}-\eps_n,
\]
that is, $\gamma_n(t)\in X_{\eps_n}$. 
We have obtained that
\[
\left\{t:I_{\eps_n}(\widetilde\gamma_n(t))>c_{\eps_n}-\frac{\eps_n}{2}\right\}
\subset\{t:\gamma_n(t)\in X_{\eps_n}\}.
\]
Now, in view of the above and of \eqref{maxequiv}, 
the asserted  equivalence~\eqref{max} follows and the claim is established.
\par
Now we recall that for every $t$ such that $\gamma_n(t)\in X_{\eps_n}$,
we have $\|\gamma_n(t)\|\le C_1$,
where $C_1$ is the constant obtained in Lemma~\ref{lem:monoest}. 
Consequently, we obtain from Lemma~\ref{lem:G} that
for every such $t$ the following estimate holds:
\begin{align}
\label{I'}
\|I_{\eps_n}'(\gamma_n(t))\|\ge\|I'(\gamma_n(t))\|-\eps_n\|G'(\gamma_n(t))\|
\ge\delta-\eps_nC_M\exp\left\{\frac{1}{8\pi}C_1^2\right\}\ge\frac{\delta}{2},
\end{align}
for sufficiently large $n$.
Moreover, for every $t$ such that $\gamma_n(t)\in X_{\eps_n}$
we also have $\Phi_n(\gamma_n(t))=1$ and $I(\gamma_n(t))\le c+\eps_n$. 
Therefore, using \eqref{exp} and \eqref{I'}
and the monotonicity property $I_{\eps_n}(\gamma_n(t))\le I(\gamma_n(t))$, we estimate
for all $\gamma_n(t)\in X_{\eps_n}$:
\begin{align*}
I_{\eps_n}(\widetilde\gamma_n(t))
\le&I_{\eps_n}(\gamma_n(t))-\sqrt{\eps_n}\Phi_n(\gamma_n(t))\|I_{\eps_n}'(\gamma_n(t))\|
+\frac{\eps_n}{2}\Phi_n(\gamma_n(t))^2\\
=&I_{\eps_n}(\gamma_n(t))-\sqrt{\eps_n}\|I_{\eps_n}'(\gamma_n(t))\|
+\frac{\eps_n}{2}\\
\le&I(\gamma_n(t))-\frac{\delta\sqrt{\eps_n}}{2}+\frac{\eps_n}{2}
\le c+\eps_n-\frac{\delta\sqrt{\eps_n}}{2}+\frac{\eps_n}{2}\\
=&c_{\eps_n}-\frac{\delta\sqrt{\eps_n}}{2}+\left(|c'|+\frac{3}{2}+\circ(1)\right)\eps_n\\
\le&c_{\eps_n}-\frac{\delta\sqrt{\eps_n}}{4}.
\end{align*}
Therefore, in view of \eqref{max}, we derive
\[
\max_{t\in[0,1]}I_{\eps_n}(\widetilde\gamma_n(t))
=\max_{\gamma_n(t)\in X_{\eps_n}}I_{\eps_n}(\widetilde\gamma_n(t))
\le c_{\eps_n}-\frac{\delta\sqrt{\eps_n}}{4}
\]
which in turn implies that
\[
c_{\eps_n}\le c_{\eps_n}-\frac{\delta\sqrt{\eps_n}}{4},
\]
a contradiction.
It follows that such a $\delta$ does not exist, and that there necessarily exists 
a bounded Palais-Smale sequence at the level $c=c(\lambda_1,\lambda_2)$.
\end{proof}
\begin{proof}[Proof of Proposition~\ref{prop:mono}]
Fix $(\lambda_1,\lambda_2)$ so that $c(\lambda_1+\lambda,\lambda_2+\lambda)$
is differentiable with respect to $\lambda$ at $\lambda=0$.
Let $(u_n)$ be the bounded Palais-Smale sequence for
$I=I_{\lambda_1,\lambda_2}$ as obtained in Lemma~\ref{lem:PS}.
We claim that $(u_n)$ converges strongly in $H^1$ to a solution $u\not\equiv0$
for problem \eqref{pb}.
Indeed, since $\|u_n\|\le C_1$, there exists $u\in E$ such that
$u_n\weak u$ weakly in $H^1$, strongly in $L^p$ for every $p\in[1,+\infty)$
and a.e.\ in $M$. 
Furthermore, $e^{u_n}\to e^u$ strongly in $L^p$ for every $p\in[1,+\infty)$
and $e^{-u_n}\to e^{-u}$ strongly in $L^p$ for every $p\in[1,+\infty)$.
Therefore, taking limits in the equation
\[
\int_M\nabla_g u_n\cdot\nabla_g\phi\,\d v_g
=\lambda_1\frac{\int_M e^{u_n}\phi\,\d v_g}{\int_M e^{u_n}\,\d v_g}
-\lambda_2\frac{\int_M e^{-u_n}\phi\,\d v_g}{\int_M e^{-u_n}\,\d v_g}
+\circ(1)
\qquad\forall\phi\in E,
\]
we see that $u$ is a solution for \eqref{pb}.
Moreover, since $\|I'(u_n)\|\to0$,
\[
\circ(1)(\|u_n-u\|)=\bra I'(u_n),u_n-u\ket=\|u_n-u\|^2+\circ(1)
\]
and therefore, $u_n\to u$ strongly in $H^1$.
By continuity and in view of Lemma~\ref{lem:c} we derive, as $n\to+\infty$
\[
I(u_n)\to I(u)=c>0
\]
and consequently $u\not\equiv0$.
\end{proof}
Let $(\lambda_1^*,\lambda_2^*)\in\Lambda$ and let $(\lambda_{1,n},\lambda_{2,n})\in\Lambda$,
$(\lambda_{1,n},\lambda_{2,n})\to(\lambda_1^*,\lambda_2^*)$ be such that problem~\eqref{pb}
admits corresponding nonzero solutions $(u_n)$ as obtained in Proposition~\ref{prop:mono}.
We note that, in view of Lemma~\ref{lem:c} with 
\[
\omega=\bar\omega=\frac{1}{2}\left(1-\frac{\lambda_1^*+\lambda_2^*}{\mu_1(M)|M|}\right)
\]
we have, as $n\to+\infty$,
\begin{equation}
\label{c}
c(\lambda_{1,n},\lambda_{2,n})\ge
\frac{\rho_{\bar\omega}^2}{4}\left(1-\frac{\lambda_1^*+\lambda_2^*}{\mu_1(M)|M|}\right)
+\circ(1).
\end{equation}
Following the notation in \cite{OS}, 
we denote by $\mathcal S_1$, $\mathcal S_2$ the blow-up sets of $(u_n)$ and $(-u_n)$,  
respectively.
Namely, we set
\begin{align*}
&\mathcal S_1=\{x\in M\ : \exists x_n\to x\ \mathrm{s.t.}\ u_n(x_n)\to+\infty\}\\
&\mathcal S_2=\{x\in M\ : \exists x_n\to x\ \mathrm{s.t.}\ u_n(x_n)\to-\infty\}.
\end{align*}
Moreover, we set
\begin{align*}
&\mu_{1,n}=\lambda_{1,n}\frac{e^{u_n}}{\int_Me^{u_n}\,\d v_g},
&\mu_{2,n}=\lambda_{2,n}\frac{e^{-u_n}}{\int_Me^{-u_n}\,\d v_g}.
\end{align*}
Since $\int_M\d\mu_{i,n}=\lambda_{i,n}$ for $i=1,2$, there exist Radon measures $\mu_i$, $i=1,2$,
such that $\mu_{i,n}\weak\mu_i$, $i=1,2$, weakly in the sense of measures.
In particular,
\begin{equation}
\label{meas}
\lambda_{i,n}=\int_M\d\mu_{i,n}\to\lambda_i^*=\int_M\d\mu_i\qquad i=1,2.
\end{equation}
At this point the following blow-up analysis from \cite{OS}
is a key step.
\begin{Proposition}[\cite{OS}]
\label{prop:os}
For the solution sequence $(u_n)$ exactly one of the following alternatives occurs:
\begin{itemize}
\item[(i)]compactness: $\mathcal S_1\cup\mathcal S_2=\emptyset$. Up to subsequences,  
$u_n\to u^*\in E$, where $u^*$ is
solution to \eqref{pb} with $\lambda_1=\lambda_1^*$ and $\lambda_2=\lambda_2^*$.
\item[(ii)]one-sided concentration: $\mathcal S_i\neq\emptyset$ and $\mathcal S_j=\emptyset$,
$j\neq i$.
In this case $\mu_i=\sum_{x_0\in\mathcal S_i}8\pi\delta_{x_0}$.
In particular, $\lambda_i^*=8\pi\,\mathrm{card}(\mathcal S_i)$.
\item[(iii)]two-sided concentration: $\mathcal S_1\neq\emptyset$ and $\mathcal S_2\neq\emptyset$.
In this case,
$\mu_i=r_i+\sum_{x_0\in\mathcal S_i}m_i(x_0)\delta_{x_0}$
with $r_i\ge0$, $r_i\in L^1(M)\cap L_{\mathrm{loc}}^\infty(M\setminus\mathcal S_i)$
and $m_i(x_0)\ge4\pi$.
Furthermore, 
\par\noindent
(iii-a)
If there exists $x_0\in\mathcal S_i\setminus S_j$, then $r_i\equiv0$
and $m_i(x_0)=8\pi$.
\par\noindent
(iii-b)
For every $x_0\in\mathcal S_1\cap\mathcal S_2$ the following relation holds:
$(m_1(x_0)-m_2(x_0))^2=8\pi(m_1(x_0)+m_2(x_0)).$
\end{itemize}
\end{Proposition} 
At this point we can complete the proof of Theorem~\ref{thm:main}.
\begin{proof}[Proof of Theorem~\ref{thm:main}]
We begin by showing that $\mathrm{card}(\mathcal S_1\cap\mathcal S_2)=\emptyset$.
Indeed, 
in view of \eqref{meas} and Proposition~\ref{prop:os}-(iii), we obtain
\[
\lambda_1^*+\lambda_2^*=\int_M\d\mu_1+\int_M\d\mu_2\ge m_1(x_0)+m_2(x_0).
\]
On the other hand, the following property is elementary:
\begin{align*}
\min\left\{x+y:x,y\ge4\pi,(x-y)^2=8\pi(x+y)\right\}
=&\min\left\{\eta:\eta\ge8\pi+|\xi|,\eta=\frac{\xi^2}{8\pi}\right\}\\
=&4(3+\sqrt5).
\end{align*}
It follows from Proposition~\ref{prop:os}-(iii-b)
that if $x_0\in\mathcal S_1\cap\mathcal S_2$, then
\begin{align*}
\lambda_1^*+\lambda_2^*\ge4(3+\sqrt5)\pi>16\pi
>\mu_1(M)|M|
\end{align*}
and therefore $(\lambda_1^*,\lambda_2^*)\not\in\Lambda$, a contradiction.
Now suppose that $\mathcal S_1=\mathcal S_1\setminus\mathcal S_2\neq\emptyset$.
Then, in view of Proposition~\ref{prop:os}-(iii-a) we obtain
$r_1\equiv0$ and $\mu_1=8\pi\sum_{x_0\in\mathcal S_1}\delta_{x_0}$ and thus \eqref{meas} implies
$\lambda_1^*=8\pi\,\mathrm{card}(\mathcal S_1)>0$ 
and again we obtain $(\lambda_1^*,\lambda_2^*)\not\in\Lambda$, a contradiction.
Similarly, we see that $\mathcal S_2=\emptyset$.
Hence, $\mathcal S_1\cup\mathcal S_2=\emptyset$. Therefore, alternative (i) 
of Proposition~\ref{prop:os} holds, and there exists a solution $u^*\in E$ to problem~\eqref{pb}
such that $u_n\to u^*$ strongly in $H^1$.
By continuity,
$I_{\lambda_{1,n},\lambda_{2,n}}(u_n)\to I_{\lambda_1^*,\lambda_2^*}(u^*)$.
In particular, in view of \eqref{c}, we have
\[
I_{\lambda_1^*,\lambda_2^*}(u^*)\ge \liminf_{n\to+\infty}c(\lambda_{1,n},\lambda_{2,n})\ge
\frac{\rho_{\bar\omega}^2}{4}\left(1-\frac{\lambda_1^*+\lambda_2^*}{\mu_1(M)|M|}\right)>0
\]
and therefore $u^*\not\equiv0$.
\end{proof}

\end{document}